\documentclass[12pt]{article}               
\usepackage{geometry}
\usepackage{cite}
\usepackage{multirow}
\usepackage{amsmath}
\usepackage{amssymb}
\usepackage[usenames]{color}
\usepackage{graphicx}          
 \usepackage{amsthm}                              
\usepackage{epsfig}    
 \usepackage{helvet}

  \newcommand{\ch}{{\rm ch\,}}
   \newcommand{\sh}{{\rm sh\,}}

\title{\large Observer Design for a Flexible Structure \\ with Distributed and Point Sensors}

\author{Alexander Zuyev and Julia Kalosha\thanks{
Institute of Applied Mathematics and Mechanics, National Academy of Sciences of Ukraine, Sloviansk
}}
\date{}

\begin{document}

\maketitle
\thispagestyle{empty}

\begin{abstract}
The paper is devoted to the observability study of a dynamic system, which describes the vibrations of an elastic beam with an attached rigid body and distributed control actions. The mathematical model is derived using Hamilton's principle in the form of the Euler--Bernoulli beam equation with hinged boundary conditions and interface condition at the point of attachment of the rigid body. It is assumed that the sensors distributed along the beam provide output information about the deformation in neighborhoods of the specified points of the beam. 
Based on the variational form of the equations of motion, the spectral problem for defining the eigenfrequencies and eigenfunctions of the beam oscillations is obtained.
Some properties of the eigenvalues ??and eigenfunctions of the spectral problem are investigated. 
Finite-dimensional approximations of the dynamic equations are constructed in the linear manifold spanned by the system of eigenfunctions. For these Galerkin approximations, observability conditions for the control system with incomplete information about the state are derived.
An algorithm for observer design with an arbitrary number of modal coordinates is proposed for the differential equation on a finite-dimensional manifold. 
Based on a quadratic Lyapunov function with respect to the coordinates of the finite-dimensional state vector, the exponential convergence of the observer dynamics is proved. 
The proposed method of constructing a dynamic observer makes it possible to estimate the full system state by the output signals characterizing the motion of particular point only. 
Numerical simulations illustrate the exponential decay of the norm of solutions of the system of ordinary differential equations that describes the observation error.
\end{abstract}

\section{Introduction}
Control theory for flexible structures is a thoroughly studied field from theoretical point of view. Controllability and stabilization problems for different models of flexible beams have been the subject of interest of many authors for several decades; see, for instance, \cite{B1978, K1992,LGM1999,LLS2012} and references therein.

In real engineering structures that consist of networks of interconnected beams, there is often a need to provide a feedback control which takes into account an output information about the state of the system. A general challenge is that most of the real settings have limitations for recording and processing output signals. Namely, only partial information about the motion of certain reference points is usually available for control purposes.
That is why the questions of observability and observer design become important for estimating the complete state of a dynamic model~\cite{L1964, L1966, B2008}.
Controllability and observability problems have been studied for strings, plates, beams and multi-link networks in \cite{L1989, BL1999, LL2000, MZ2002, Z2005,KGBS2008,GZR217, MTG2020}.

An analytical approach for designing a Luenberger observer for a finite-dimensional linear system is based on constructing an auxiliary matrix of gain parameters such that the  error dynamics is asymptotically stable.
Then the above gain parameters can be defined by checking the  Routh--Hurwitz or the Li{\'e}nard--Chipart conditions for the corresponding characteristic polynomial.
For mathematical models of flexible structures with multiple degrees of freedom, this approach requires the analysis of high order
polynomial inequalities with parameter dependent coefficients. The implementation of this scheme may be a very difficult technical task.
An alternative approach for observer design is based on Lyapunov's direct method.
 The idea is to take a suitable quadratic form as a Lyapunov function candidate and define the gain parameters
 in such a way that the observation error dynamics has stable trivial equilibrium.
Moreover, it is possible to prove that the trivial solution is asymptotically stable with the help of the Barbashin--Krasovskii theorem (or LaSalle's invariance principle).
This technique has been applied to observer design and observed-based stabilization in~\cite{ZS2005, ZS2006, ZS2007} for mathematical models of flexible manipulators.

In the present paper, we study the mathematical model that describes a flexible beam attached to a rigid body (shaker). The beam is simply supported at the ends by two hinges and is actuated by distributed piezoelectric controllers and the shaker. In previous studies, the stabilization problem with a state feedback was investigated for this infinite-dimensional system. The control design that ensures strong asymptotic stability of the trivial equilibrium was proposed in~\cite{KZ2021}. Now, we widen our study to the system with output and question ourselves whether it is possible to fully reconstruct the state vector having at our disposal only output signals on a subregion of the domain. In this work, the observability problem is investigated for a finite-dimensional approximation of the dynamics. An observer design that guarantees the exponential decay of the observation error is proposed.

The rest of this paper is organized as follows. The mechanical system under con\-sideration is described in section~2, and its equations of motion are derived from the variational principle in section~3. Some properties of eigenvalues of the associated spectral problem are discussed in sections~4.
A finite-dimensional approximation of the dynamical equations based on Galerkin's method is proposed in section~5.
The main results (observability conditions together with observer design scheme) are presented in section~6.
Finally, we illustrate the convergence of the proposed observer with some numerical simulations in section~7.

\section{Modeling of a flexible beam with a shaker}

Consider a flexible beam supported at both ends and attached to the shaker. Let $l$ be the length of the beam, and let $w(x,t)$ denote the deflection of the centerline of the beam at a point $x\in [0,l]$ and time $t$ (here $x=0$ stands for the upper end and $x=l$ for the lower one, correspondingly). We assume that the shaker is attached at the point $x=l_0$, and $m$ is the mass of the moving part of the shaker. Then the kinetic energy of the system considered is
\begin{equation}
K=\frac12\int_0^l \rho(x) {\dot w}^2(x,t)\,dx + \frac{m}2 \dot w^2(l_0,t),
\label{kinetic_K}
\end{equation}
where $\rho(x)$ is the mass per unit length of the beam. We use dots to denote derivatives with respect to $t$, and primes to denote derivatives with respect to $x$. According to the Euler--Bernoulli beam model, the potential energy can be written as follows:
\begin{equation}
U=\frac12\int_0^l E(x)I(x) \left(w''(x,t)\right)^2\,dx + \frac{\varkappa}2 w^2(l_0,t).
\label{potential_U}
\end{equation}
Here $E(x)$ and $I(x)$ are the Young modulus and the moment of inertia of the cross section area, respectively, and $\varkappa$ is the stiffness coefficient of the shaker.
In the sequel, we assume that damping may be neglected.

\section{Variational form of the dynamical equations}
Suppose that $w(x,t)$ defines the motion of the system considered for $t\in [0,\tau]$, $\tau>0$. We have the following geometric boundary conditions:
$$
w(0,t)=w(l,t)=0.
$$
By using Hamilton's principle, we get
\begin{equation}
\delta \int_0^\tau (K-U)\,dt + \int_0^\tau \left\{F\delta w(l_0,t)+\sum_{j=1}^k M_j\int_0^l\psi_j(x)\delta w''(x,t)dx\right\}dt=0,
\label{Hamilton_beam}
\end{equation}
for each admissible variation $\delta w(x,t)$ of class $C^2\left([0,l]\times [0,\tau]\right)$ satisfying the boundary conditions
$$
\delta w|_{t=0}=\delta w|_{t=\tau}=\delta w|_{x=0}=\delta w|_{x=l}=0.
$$
Here $F$ is the control force implemented by the shaker. In formula~\eqref{Hamilton_beam}, we assume also that there are $k$ piezoelectric actuators attached to the beam. We describe the action of the $j$-th actuator by its density of the torque $M_j$ and the shape function $\psi_j(x)$. These characteristics may be computed by using results of the paper \cite{AT}. We assume that
$$
{\rm supp}\, \psi_j \cap \{0,l_0,l\} = \emptyset \quad \text{for each}\; j=1,2,\dots,k.
$$

By performing the integration by parts with respect to $t$ in~\eqref{Hamilton_beam}, we obtain
\begin{equation}\label{VF}
\int_0^l \left(\rho \ddot w \delta w+ \left(EI w'' -\sum_{j=1}^kM_j\psi_j\right)\delta w''\right)dx+\left.\left(m\ddot w+\varkappa w - F\right)\delta w\right|_{x=l_0}=0.
\end{equation}

Then another integration by parts with respect to $x$ yields the following partial differential equation:
\begin{equation}
\label{E_B}
\rho \ddot w(x,t) + \left(EI w''(x,t)\right)'' = \sum_{j=1}^k \psi_j''(x)M_j,\quad x\in (0,l) \setminus\{l_0\},
\end{equation}
with boundary and interface conditions
\begin{equation}\label{BC}
\begin{gathered}
w|_{x=0}=w|_{x=l}=0,\;w''|_{x=0}=w''|_{x=l}=0,\\
\left.\left(m\ddot w+\varkappa w\right)\right|_{x=l_0} = \left.(EIw'')'\right|_{x=l_0-0} - \left.(EIw'')'\right|_{x=l_0+0}+ F,
\end{gathered}
\end{equation}
for $w(x,t)$ of class $C^2\left([0,l]\times [0,\tau]\right)$ such that
the derivatives $w'''(x,t)$ and $w''''(x,t)$ exist for all $x\in[0,l]\setminus\{l_0\}$ and $t\in[0,\tau]$.

\section{Spectral problem}
To analyse the spectrum of oscillations of the mechanical system considered, we first study the control system~\eqref{E_B}--\eqref{BC} with ${M_1=\ldots=M_k=0}$ and $F=0$. For simplicity, we assume that the beam is homogeneous, i.e. $\rho$, $E$, and $I$ are treated as constants in the sequel.

By substituting $w(x,t)=W(x)q(t)$ into the homogeneous part of~\eqref{E_B}--\eqref{BC}, we get $\ddot q(t)=-\lambda q(t)$ together with the following spectral problem:
\begin{equation}\label{spectral}
\begin{gathered}
\frac{d^4}{dx^4}W(x) = \lambda \frac{\rho}{EI}W(x),\quad x\in (0,l)\setminus \{l_0\},\\
W(0)=W(l)=0,\; W''(0)=W''(l)=0,\;W\in C^2[0,l],\\
W'''(l_0-0)-W'''(l_0+0)=\frac{\varkappa-\lambda m}{EI}W(l_0).
\end{gathered}
\end{equation}
The eigenfunctions of the above problem are
$$
W(x)=\left\{\begin{array}{ll}C_1 \sin(\mu x)+C_2\sinh(\mu x),& x\le l_0,\\
B_1 \sin(\mu (x-l))+B_2\sinh(\mu (x-l)),& x> l_0,
\end{array}\right.
$$
where $\mu = \left(\frac{\lambda\rho}{EI}\right)^{1/4}$, and the constants $C_1, C_2, B_1, B_2$ satisfy the following linear algebraic system:
{\small
\begin{equation}
\left(\begin{array}{cccc}\sin \mu l_0 & -\sh \mu l_0 & -\sin \mu(l_0-l) & \sh\mu(l_0-l) \\
\sin \mu l_0 & \sh \mu l_0 & -\sin \mu(l_0-l) & -\sh\mu(l_0-l)\\
\cos \mu l_0 & \ch \mu l_0 & -\cos \mu(l_0-l) & -\ch\mu(l_0-l)\\
\cos \mu l_0 + \frac{\varkappa-m\lambda}{EI\mu^3}\sin \mu l_0&  \frac{\varkappa-m\lambda}{EI\mu^3}\sh \mu l_0 -\ch \mu l_0& -\cos \mu(l_0-l) & \ch\mu(l_0-l)\end{array}\right)\left(\begin{array}{c}C_1 \\ C_2 \\ B_1 \\ B_2\end{array}\right) =0.
\label{algsys_shaker}
\end{equation}}
Thus, the eigenvalues $\lambda$ of the spectral problem~\eqref{spectral} are obtained by solving the equation
\begin{equation}
\Delta(\lambda) = 0,
\label{chareq_shaker}
\end{equation}
where $\Delta(\lambda)$ is the determinant of the $4\times 4$-matrix of~\eqref{algsys_shaker}.
A rigorous investigation of the characteristic equation~\eqref{chareq_shaker} with regard to the asymptotic behavior of its roots is carried out in~\cite{KZB2021}.

Let $\lambda_1$, $\lambda_2$, $\dots$ be eigenvalues of~\eqref{spectral}, and let $W_1(x)$ ,$W_2(x)$, $\dots$  be the corresponding eigenfunctions. Let us define
\begin{equation}
\left<W_i,W_j\right>_H:=\int_0^l \rho W_i(x)W_j(x)\,dx + m W_i(l_0) W_j(l_0).
\label{bilinearform}
\end{equation}
We have the following lemma.

{\bf Lemma~1.}
{\em
If $\lambda_i\neq \lambda_j$ then $\left<W_i,W_j\right>_H=0$.}

{\bf Proof.}
The integration by parts of the expression $\int_0^lW_i''(x)W_j''(x)dx$ and exploiting the interface conditions $W_j'''(l_0+0)-W_j'''(l_0-0)=\frac{\lambda m - \varkappa}{EI}W_j(l_0)$, $j=1,\dots,N$, gives
$$\int_0^l W_j(x)\frac{d^4}{dx^4}W_i(x)dx+\frac{\lambda_im-\varkappa}{EI}W_i(l_0)W_j(l_0)$$
on the one hand, and
$$\int_0^l W_i(x)\frac{d^4}{dx^4}W_j(x)dx+\frac{\lambda_jm-\varkappa}{EI}W_i(l_0)W_j(l_0)$$
on the other hand. Calculating the difference between these two expressions, we obtain
$$\int_0^l\left(W_j(x)\frac{d^4}{dx^4}W_i(x)-W_i(x)\frac{d^4}{dx^4}W_j(x)\right)dx+\frac m{EI}(\lambda_i-\lambda_j)W_i(l_0)W_j(l_0)=0.$$
As both $W_i(x)$ and $W_j(x)$ satisfy~\eqref{spectral},
$$(\lambda_i-\lambda_j)\left<W_i,W_j\right>_H=0,$$
which leads to $\left<W_i,W_j\right>_H=0$ for $\lambda_i\neq\lambda_j$.
$\square$

{\bf Lemma~2.}
{\em
Each eigenvalue of the spectral problem~\eqref{spectral} is a non-negative real number.}

{\bf Proof.}
Consider the differential operator ${\cal L}:D({\cal L})\to \tilde H = L^2[0,l]\times {\mathbb R}$ acting as
$$
\zeta=\left(
                   \begin{array}{c}
                     u(x) \\
                     p \\
                   \end{array}
                 \right)\mapsto {\cal L}\zeta=\left(
                                            \begin{array}{c}
                                              \frac{EI}\rho u''''(x) \\
                                              \frac\varkappa mp-\frac{EI}m\left(u'''(l_0-0)-u'''(l_0+0)\right) \\
                                            \end{array}
                                          \right)$$
with the domain
$$D({\cal L})=\left\{
\begin{aligned}
&\zeta\in H^2[0,l]\cup H^4(0,l_0)\cup H^4(l_0,l)\times{\mathbb R}:\\
&u(0)=u(l)=0,\;u''(0)=u''(l)=0,\\
&u''(l_0-0)=u''(l_0+0),\;p=u(l_0)
\end{aligned}
\right\}\subset \tilde H.$$
Then the spectral problem~\eqref{spectral} can be rewritten as
$${\cal L}\zeta=\lambda\zeta.$$
Let us rewrite the bilinear form~\eqref{bilinearform} in terms of elements of $\tilde H$:
$$\left<\zeta_i,\zeta_j\right>_{\tilde H}=\left<
\left(
  \begin{array}{c}
    u_i \\
    p_i \\
  \end{array}
\right),               \left(
                         \begin{array}{c}
                           u_j \\
                           p_j \\
                         \end{array}
                       \right)\right>_{\tilde H}:=\int\limits_0^l\rho u_i(x)u_j(x)dx+mp_ip_j.$$
A direct computation leads to
$$\left<{\cal L}\zeta_i,\zeta_j\right>_{\tilde H}=\int\limits_0^lEIu_i''''(x)u_j(x)dx+\varkappa p_ip_j-EI\left(u_i'''(l_0-0)-u_i'''(l_0+0)\right)p_j$$
and
$$\left<\zeta_i,{\cal L}\zeta_j\right>_{\tilde H}=\int\limits_0^lEIu_i(x)u_j''''(x)dx+\varkappa p_ip_j-EI\left(u_j'''(l_0-0)-u_j'''(l_0+0)\right)p_i.$$
The integration by parts gives
$\left<{\cal L}\zeta_i,\zeta_j\right>_{\tilde H}=\left<\zeta_i,{\cal L}\zeta_j\right>_{\tilde H}$.
Besides, it is easily verifiable that
$\left<\zeta,{\cal L}\zeta\right>_{\tilde H}=\int\limits_0^lEI(u''(x))^2dx+\varkappa p^2\geq0$ for all $\zeta\in D({\cal L})$.
So, ${\cal L}$ is positive self-adjoint operator with respect to the bilinear form $\left<\cdot,\cdot\right>_{\tilde H}$; thus the eigenvalues of ${\cal L}$ are real and non-negative.
$\square$


\section{Galerkin's method}
Let us consider a finite set of eigenvalues $\lambda_1$, $\lambda_2$, ... , $\lambda_N$ and the corresponding eigenfunctions
$W_1(x)$, $W_2(x)$, ..., $W_N(x)$ of the spectral problem~\eqref{spectral}.
Then we substitute
$$
w(x,t) = \sum_{j=1}^N W_j(x)q_j(t)
$$
into the variational form of equations of motion~\eqref{VF}, assuming that~\eqref{VF} holds for each
$$
\delta w(\cdot,t)\in {\rm span}\left\{W_1,W_2,\dots,W_N\right\},\quad t\ge 0.
$$
As a result, we get the Galerkin system with respect to $q_1(t),q_2(t),\dots,q_N(t)$:
\begin{equation}
\label{Galerkin1}
\ddot q_j + \lambda_j q_j = \frac{W_j(l_0)}{\|W_j\|^2_H}u_0+\sum_{i=1}^k \frac{\int_0^l\psi_i(x) W_j''(x)\,dx}{\|W_j\|^2_H}u_i,\;\; j=1,2,\dots,N,
\end{equation}
where $u_0=F$, $u_i=M_i$, $i=1,\dots,k$.

System~\eqref{Galerkin1} can be equivalently written as
\begin{equation}
\dot z= Az+Bu,\quad z=(q_1,p_1,\dots,q_N,p_N)^T,\;u=(u_0,u_1,\dots,u_k)^T,
\label{Galerkin3}
\end{equation}
where $A = {\rm diag}(A_1,A_2,\dots,A_N)$,
$$A_j = \left(\begin{array}{cc}0 & 1 \\ -\lambda_j & 0\end{array}\right),\;
B = \left(\begin{array}{cccc}0 & 0 & \ldots & 0 \\
b_{10} & b_{11} & \ldots & b_{1k} \\
\vdots  & \vdots  & \ddots & \vdots \\
0 & 0 & \ldots & 0 \\
b_{N0} & b_{N1} & \ldots & b_{Nk}\end{array}\right),
$$
\begin{equation}
b_{j0}=\frac{W_j(l_0)}{\|W_j\|^2_H},\; b_{ji}=\frac{\int_0^l\psi_i(x) W_j''(x)\,dx}{\|W_j\|^2_H},\; i=1,2,\dots,k.
\label{b_coeff}
\end{equation}

\section{Luenberger observer design}
Assume that the state $w(x,t)$ of the flexible beam can measured by $r$ strain gauges located at points $x=l_1,\dots,l_r$ and, moreover, the shaker displacement $w(l_0,t)$ is measured as function of time.
 The current, generated by the $s$-th gauge, is proportional to $\dot w''(l_s,t)$.
By assuming that the measuring device implements a current integrator, we treat the signals
\begin{equation}
y_0(t) = w(l_0,t),\;
y_s(t) = w''(l_s,t),\quad s=1,\dots,r,
\label{output_shaker}
\end{equation}
as the output for the mathematical model~\eqref{E_B}--\eqref{BC}. Then, as the representation
$$
w(x,t) = \sum_{j=1}^N W_j(x)q_j(t)
$$
is used to derive Galerkin's system~\eqref{Galerkin3}, we express the output~\eqref{output_shaker} in terms of functions $q_1(t), \dots, q_N(t)$ as follows:
$$
y_0(t)\approx \sum_{j=1}^N W_j(l_0)q_j(t),\; y_s(t)\approx \sum_{j=1}^N W_j''(l_s)q_j(t),\quad s=1,\dots,r.
$$

Thus, we assume that the output signal $y(t)=Cz(t)\in {\mathbb R}^{r+1}$ is available for the finite-dimensional control system~\eqref{Galerkin3},
where
\begin{equation}
y = C z,\,
C=\left(\begin{array}{lllllll}
c_{01} & 0 &  \ldots & c_{0N} & 0 \\
c_{11} & 0 &  \ldots & c_{1N} & 0\\
\vdots & \vdots  & \ddots & \vdots & \vdots\\
c_{r1} & 0 &  \ldots & c_{rN} & 0\end{array}\right),\;c_{0j} = W_j(l_0),\, c_{sj}=W_j''(l_s),\,s=\overline{1,r}.
\label{outp}
\end{equation}

{\bf Lemma~3.}
{\em
Let the eigenvalues $\lambda_1$, $\lambda_2$, ..., $\lambda_N$ be distinct
and, for each $j=\overline{1,N}$, there is an $s\in\{0,...,r\}$ such that $c_{sj}\neq0$.
Then system~\eqref{Galerkin3} with output~\eqref{outp} is observable.
}

{\bf Proof.}
Let us prove that ${\rm rank}\,H=2N$, where the matrix
$$H=\left(
\begin{array}{l}
C \\
CA \\
\vdots \\
CA^{2N-1}
\end{array}\right)
$$
by a suitable permutation of its rows and columns can be transformed to the following block form:
$$\left(
    \begin{array}{cc}
      H_0 & 0 \\
      \vdots & \vdots \\
      H_r & 0 \\\hline
      0 & H_0 \\
      \vdots & \vdots \\
      0 & H_r \\
    \end{array}
  \right),\quad
H_s=\left(
      \begin{array}{ccc}
        c_{s1} & \ldots & c_{sN} \\
        -\lambda_1c_{s1} & \ldots & -\lambda_Nc_{sN} \\
        \vdots & \ddots & \vdots \\
        -\lambda_1^{N-1}c_{s1} & \ldots & -\lambda_N^{N-1}c_{sN} \\
      \end{array}
    \right),\quad s=0,\dots,r.$$
Each $H_s$ is a Vandermonde-like matrix whose determinant can be calculated as
$$\det H_s=\prod\limits_{i=1}^Nc_{si}\prod\limits_{1\leq j<n\leq N}(\lambda_n-\lambda_j),\quad s=0,\dots,r.$$
Under the assumptions of Lemma~3,
the rectangular block
$\left(\begin{array}{cc}
          H_0 \\
          \vdots \\
          H_r \\
       \end{array}\right)$
has a nonzero minor of order $N$. So, ${\rm rank}H=2N$ and system~\eqref{Galerkin3},~\eqref{outp} is observable by Kalman's observability criterion.
$\square$

The following result describes an explicit scheme for constructing the Luenberger observer for system~\eqref{Galerkin3},~\eqref{outp}.

{\bf Theorem~1.}
{\em
Let system~\eqref{Galerkin3},~\eqref{outp} be observable and $\lambda_j>0$ for all $j=\overline{1,N}$. Then, for any $z(0)\in {\mathbb R}^{2N}$, $\bar z(0) \in {\mathbb R}^{2N}$, and any $u\in L^\infty[0,\infty)$, the corresponding solutions $z(t)$ and $\bar z(t)$ of system~\eqref{Galerkin3} and
\begin{equation}\label{observer}
\dot {\bar z}(t) = (A-FC)\bar z(t) + B u(t) + F y(t)
\end{equation}
satisfy the property
\begin{equation}
\|z(t)-\bar z(t)\| \to 0\quad \text{as}\;
 t\to+\infty.
\label{errconv}
\end{equation}
Here
$$F=\left(
           \begin{array}{ccc}
             f_{10} & \ldots & f_{1r} \\
             0 & \ldots & 0 \\
             \vdots & \cdots & \vdots \\
             f_{N0} & \ldots & f_{Nr} \\
             0 & \ldots & 0 \\
           \end{array}
         \right),\quad f_{js}=\gamma_s\frac{c_{sj}}{\lambda_j\|W_j\|^2},\quad j=\overline{1,N},\;s=\overline{0,r},$$
and $\gamma_s>0$ are arbitrary constants.
}

{\bf Proof.}
Consider the observation error $e(t)=z(t)-\bar z(t)=(\Delta_1,\delta_1,\ldots,\Delta_N,\delta_N)^T$. Then subtracting~\eqref{observer} from~\eqref{Galerkin3} yields the following dynamics:
\begin{equation}\label{error}
\dot e=(A-FC)e.
\end{equation}

Let us consider the following positive definite quadratic form:
\begin{equation*}
2W(e)=\sum\limits_{j=1}^N\|W_j\|^2(\delta_j^2+\lambda_j\Delta_j^2),
\end{equation*}
and compute its time derivative along the trajectories of system~\eqref{error}:
$$\dot W(e)=-\sum\limits_{s=1}^r\sum\limits_{j=1}^N \|W_j\|^2\lambda_j\Delta_jf_{js} \sum\limits_{i=1}^N C_{si}\Delta_i = -\sum\limits_{s=1}^r\gamma_s\left(\sum\limits_{j=1}^N C_{sj}\Delta_j\right)^2\leq0.$$
According to Lyapunov's theorem, the trivial solution of~\eqref{error} is stable.

Our goal is to check whether the system~\eqref{error} admits nontrivial trajectories on the set ${\cal K}=\{e\in\mathbb R^{2N}:\dot W(e)=0\}$.
As far as $\gamma_s>0$ for all $ s=0,\dots,r$, the identity $\dot W(e(t))\equiv0$ implies that $\sum\limits_{j=1}^N C_{sj}\Delta_j(t)\equiv0$ or, in the matrix form,
\begin{equation}\label{attr1}
Ce(t)\equiv0,
\end{equation}
which in turn implies that the time derivatives $\frac{d^n}{dt^n}(Ce(t))$ vanish for all $n=0,1,2,...$.
In view of~\eqref{attr1}, system~\eqref{error} reduces to
\begin{equation}\label{attr2}
\dot e=Ae
\end{equation}
on the set ${\cal K}$. Computing the derivatives of $Ce(t)$ up to $(2N-1)$-st along the trajectories of~\eqref{attr2}, one can see that
$$\frac{d^n}{dt^n} Ce(t)=CA^ne(t),\quad n=0,\ldots,2N-1.$$
So, the components of $e(t)$ satisfy the following algebraic system:
\begin{equation}\label{attr_alg}
CA^ne(t)=0,\quad n=0,\ldots,2N-1.
\end{equation}
As system~\eqref{Galerkin3},~\eqref{outp} is observable, system~\eqref{attr_alg} admits only the trivial solution $e(t)\equiv0$ because of Kalman's observability rank condition. Thus, the trivial solution of~\eqref{error} is asymptotically stable by the Barbashin--Krasovskii theorem.
$\square$

{\em Remark.} As the error dynamics~\eqref{error} is linear in the considered case, the obtained asymptotic stability result
guarantees in fact the exponential convergence in~\eqref{errconv}.
Note that the assertion of Theorem~1 can be reformulated in terms of partial asymptotic stability of the extended system~\eqref{Galerkin3},~\eqref{outp},~\eqref{observer} with respect to the variables $z-\bar z$.
This type of stability conditions can be formally analyzed with extensions of the Barbashin--Krasovskii--LaSalle results for partial stabilization of ordinary differential equations (see, e.g.,~\cite{Z2001}) and dynamical systems in abstract spaces~\cite{Z2006}.

\section{Numerical simulations}
The goal of this section is to illustrate the behavior of solutions of control system~\eqref{Galerkin3},~\eqref{outp}
and its dynamic observer~\eqref{observer} by numerical simulations for the model with $N=6$ modes of oscillations.
All computations are carried out for the following realistic mechanical parameters (cf.~\cite{KZB2021}):
$$\begin{aligned}
l=1.875\,\text{m},\quad l_0=1.4\,\text{m},\quad E=7.1\cdot10^{10}\,\text{Pa},\quad I=1.6875\cdot10^{-10}\,\text{m}^4,\\ \rho=2660\cdot2.25\cdot10^{-4}\,\text{kg}/\text{m},\quad m=0.045\,\text{kg},\quad \varkappa=2630\,\text{N/m}.
\end{aligned}$$
The model is driven by the periodic force $F=\sin4t$ and the beam displacement at $x=l_0$ is taken as the output,
so that the piezo actuators are switched off the data from piezo sensors is neglected.
For numerical simulations, we choose the observer gain parameter $\gamma_1=3$ and the initial conditions $\bar z(0)=0$, $q_i(0)=0.1$, $p_i(0)=0.1$, $i=1,\dots,N$.
The overall weighted error $\|e(t)\|^2=\frac\varkappa2\sum\limits_{j=1}^N(\Delta_j(t))^2+\frac m2\sum\limits_{j=1}^N(\delta_j(t))^2$ is presented in Fig.~1, while high frequency components of the observation error are depicted in Figs.~2 and~3.
\begin{figure}[ht]
\centering
\includegraphics[scale=0.3]{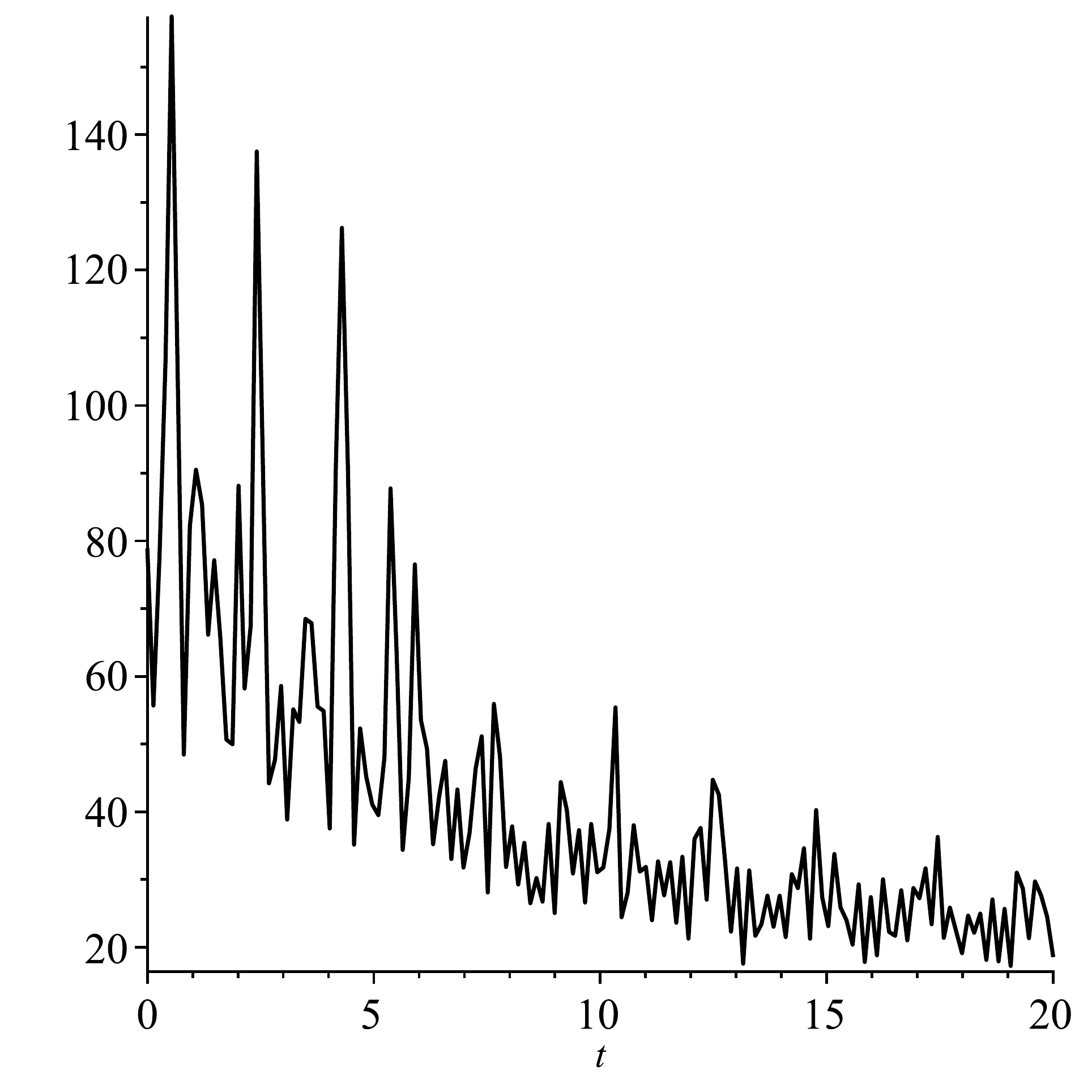}
\\
{\small Fig.~1. The time plot of $\|e(t)\|^2$.}
\end{figure}
\begin{figure}[ht]
\parbox{.45\linewidth}{\centering
\includegraphics[scale=0.275, clip]{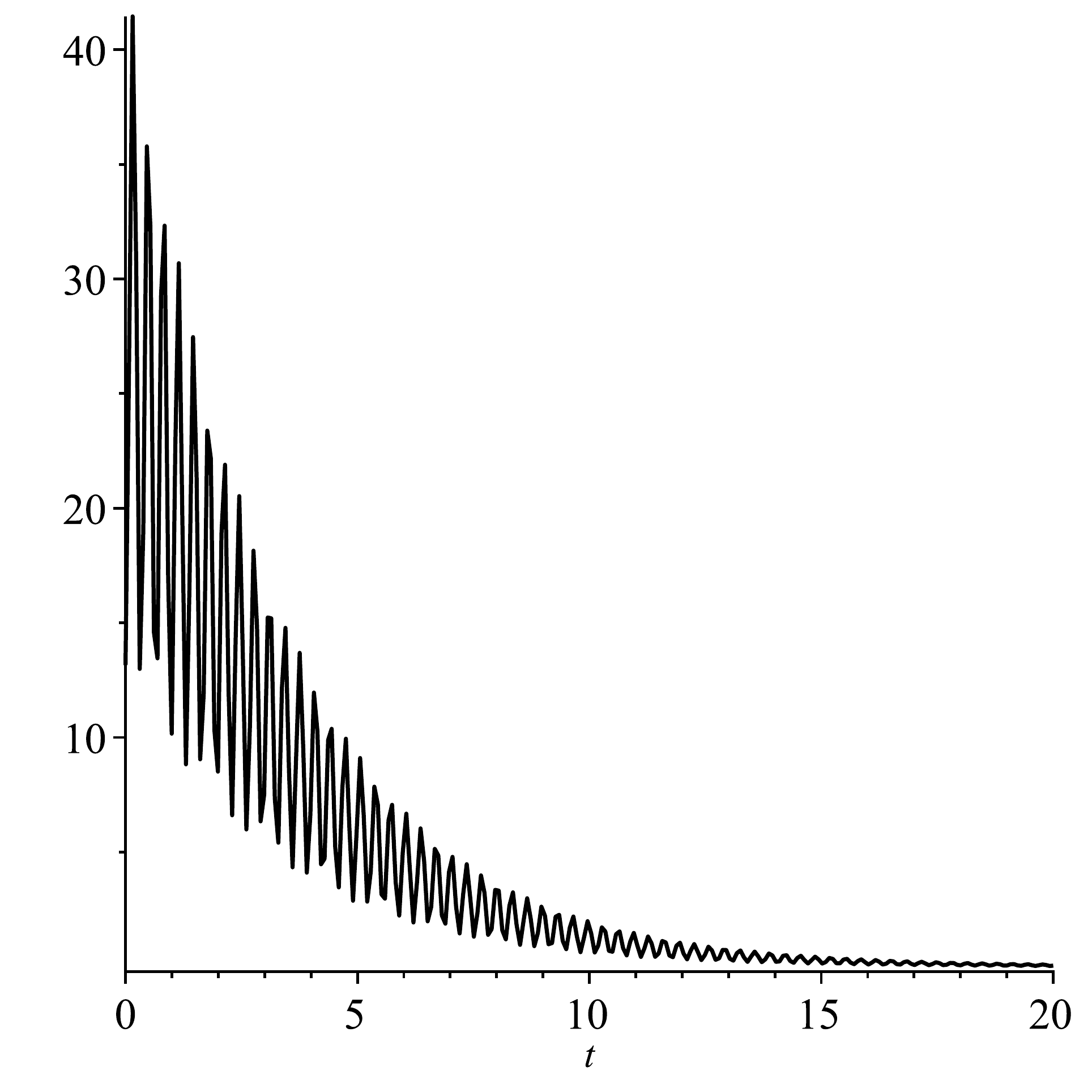}
\\
\small Fig.~2. The time plot of {$\frac\varkappa2(\Delta_N(t))^2+\frac m2(\delta_N(t))^2$.}}
\parbox{.2\linewidth}{}
\parbox{.45\linewidth}{\centering
\includegraphics[scale=0.275, clip]{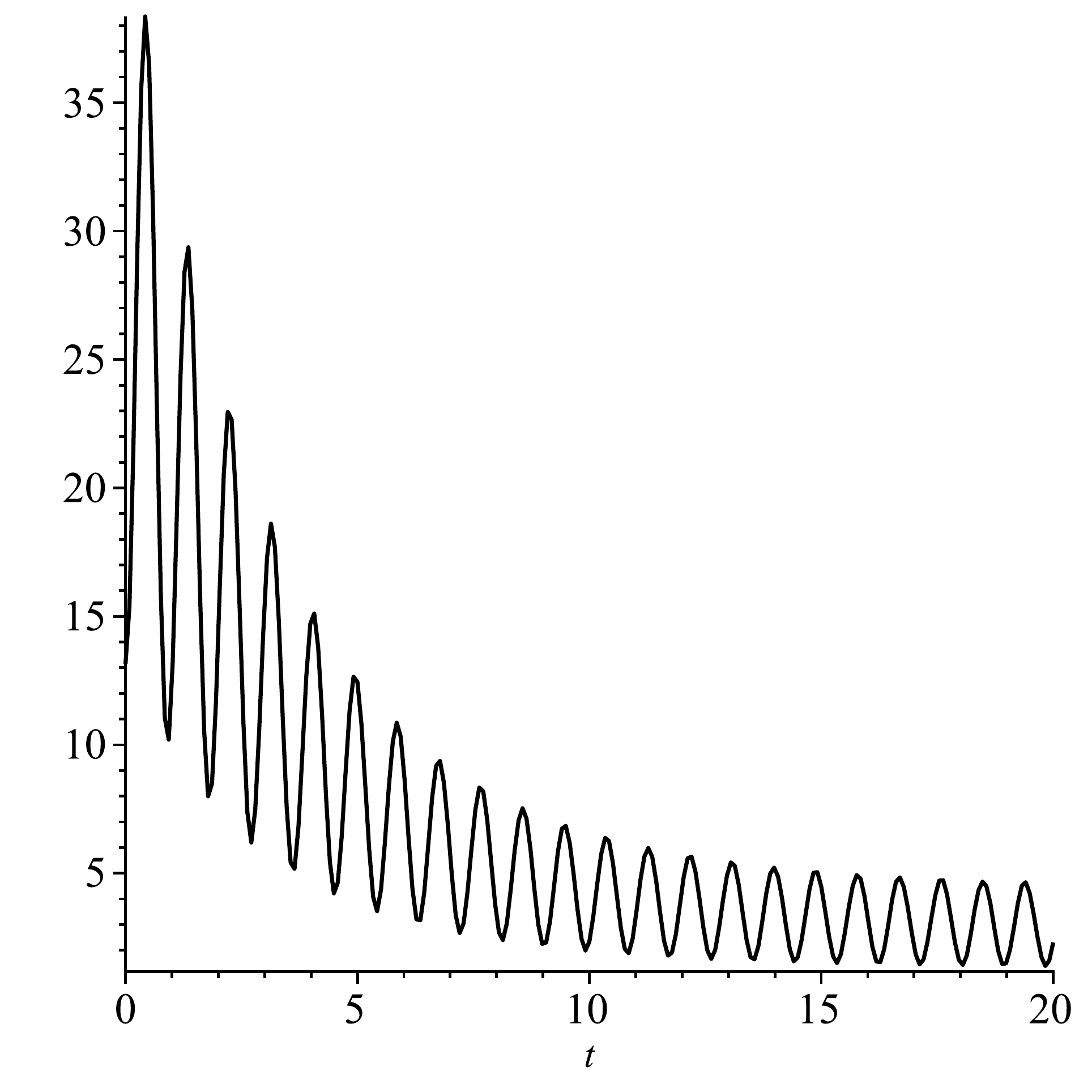}
\\
{\small Fig.~3. The time plot of $\frac\varkappa2(q_N(t)-\bar q_N(t))^2+\frac m2(p_N(t)-\bar p_N(t))^2$.}
}
\end{figure}

The above time plots confirm the convergence of the proposed observer for the considered multidimensional model.

\section{Conclusions and future work}
In this paper, we have proposed an explicit analytic approach for constructing the Luenberger observer for
a class of finite-dimensional models of flexible structures.
This observer design allows estimating the full system state for an arbitrary dimension of the output vector.
We do not solve here  the output stabilization problem leaving this task for future investigations.
In the previous work~\cite{KZ2021}, the stabilization problem was solved for the beam-shaker model
by a feedback law depending on the infinite-dimensional state vector.
So, another direction of further study is related to the observer-based stabilization of this class of distributed parameter flexible structures.

\end{document}